\documentclass[a4paper,12pt]{amsart}
\usepackage[T1]{fontenc}    
\usepackage[utf8]{inputenc}
\usepackage{enumerate}
\usepackage{amssymb}
\usepackage{fullpage}
\usepackage{tikz, subfigure}
\usetikzlibrary{patterns}

\def\RR{\mathbb R}

\def\zz{\mathbb Z}

\def\nn{\mathbb N}

\newcommand{\set}[1]{\left\lbrace #1\right\rbrace}

\newcommand{\remove}[1]{ }
\newcommand{\qtq}[1]{\quad\text{#1}\quad}

\newtheorem{theorem}{Theorem}[section]
\newtheorem{lemma}[theorem]{Lemma}

\newtheorem*{theorem-1}{Theorem 1.1}

\theoremstyle{definition}

\theoremstyle{remark}
\newtheorem*{remark}{Remark}
\remove{}

\numberwithin{equation}{section}
\begin{document}
\title[]{  Intersections of Cantor sets with hyperbolas and  continuous images}
\author{Yi Cai}
\address[Y. Cai]{School of Sciences, Shanghai Institute of Technology, Shanghai 201418,
People's Republic of China}
\email{caiyi@sit.edu.cn}

\author{Xiu Chen}
\address[X. Chen]{School of Mathematics and Physics, School of Cryptology, Nanjing Institute of Technology, Nanjing 211167, People's Republic of China}
\email{chenxiu1216@163.com}

\author{ Lipeng Wang}
\address[L. Wang]{School of Science, Jiangsu Ocean University, Lianyungang 222005, People's Republic of China}
\email{lipengwang@jou.edu.cn}

\begin{abstract}
Given $\lambda\in (0,1/2)$, let
\begin{equation*}
C_\lambda=\set{(1-\lambda)\sum_{i=1}^\infty d_i\lambda^{i-1}:d_i\in\set{0,1}}
\end{equation*}
be the middle Cantor sets with convex hull $[0, 1]$. We are interested in the set $S_t=\set{(x,y)\in C_\lambda\times C_\lambda: xy=t}$, where $t\in[0,1]$. Since the cases where $t=0$ or $t=1$ are trivial, we assume that $t\in(0,1)$ in what follows. We show that there exists a $\lambda_0=0.4302$ such that for all $\lambda$ satisfying $\lambda_0 \le \lambda < 1/2$, the set $S_t$ has the cardinality of the continuum for every $t \in (0,1)$. Besides, we further investigate the continuous image of $C_\lambda\times C_\lambda$, that is,  for any given $2\le k\in \nn$, we give a sufficient condition for set $\set{x^ky:x,y\in C_\lambda}$ to be the interval $[0,1]$. Our observations reveal that the behavior exhibited by the image of the function $f_k(x,y)=x^ky$ is complex and depends on the parameters $k$ and $\lambda$.

\end{abstract}
\maketitle
\section{Introduction}
Let $\lambda\in (0,1/2)$, and let $C_\lambda$ be the unique nonempty compact set in $\RR$ that satisfies
\begin{equation*}
  C_\lambda=\varphi_0(C_\lambda)\cup \varphi_1(C_\lambda),
\end{equation*}
where $\varphi_i(x)=\lambda x+i(1-\lambda),i=0,1$(cf. \cite{KF}). In this case, $C_\lambda$ is the classical middle Cantor set with convex hull $[0, 1]$.

The intersections of fractal sets with curves have been studied by many scholars. Zhang, Jiang and Li \cite{ZJL} studied the empty intersection of the lines with $C_\lambda\times C_\lambda$. Cai and Yang \cite{CC} extended the consideration to the empty intersection between non-self-similar sets and lines. Du, Jiang and Yao \cite{DJY} proved that the intersection of $C_{1/3}\times C_{1/3}$ with unit circle contains more than ten million elements.

Recently, Jiang et al. \cite{JKLW} considered the intersection of $C_\lambda\times C_\lambda$ with unit circle. They showed a surprising result, namely, if $0.407493\le \lambda<\frac{1}{2}$ the intersection is of cardinality continuum.  Motivated by their research, it is natural to investigate the behavior of the intersections of $C_\lambda\times C_\lambda$ with a family of curves.

The first part of this paper is devoted to the case of intersections between $C_\lambda\times C_\lambda$ and a family of hyperbolas. More precisely, let $t\in(0,1)$ and
\begin{equation*}
 S_t=\set{(x,y)\in C_\lambda\times C_\lambda:xy=t}.
\end{equation*}
So $S_t$ is the intersection of $C_\lambda\times C_\lambda$ with hyperbola $y=\frac{t}{x}$. We address the following question: does there exist an $\lambda$ such that $S_t$ is uncountable for any $t\in(0,1)$? Our first result answers this question affirmatively.

\begin{theorem}\label{theorem-1} If $0.4302\le \lambda<1/2$, then $S_t$ has the cardinality of the continuum for any $0<t<1$.
\end{theorem}

In the latter part, we focus on the continuous image of $C_\lambda\times C_\lambda$. This is closely related to the arithmetic and structure of Cantor sets, which have attracted substantial research attention, e.g.  \cite{BK,DKY,PN,YK,TJ,ZJ,ZLL}. Given $ k\in \nn=\set{1,2,3,\dots}$, let $f_k(x,y)=x^ky$. Define
\begin{equation*}
  f_k(C_\lambda,C_\lambda)=\set{f_k(x,y):x,y\in C_\lambda}.
\end{equation*}
The interest in this set stems from several pioneering works \cite{ART,JK,LJ}. Under some checkable conditions on the partial derivatives of continuous $g(x, y)$, Jiang \cite{JK} proved that
\begin{equation*}
  g(K_1,K_2)=\left[\min_{(x,y)\in K_1\times K_2} g(x,y),\max_{(x,y)\in K_1\times K_2} g(x,y)\right],
\end{equation*}
where $K_1,K_2$  be two Cantor sets with convex hull $[0,1]$. As an application of this result, he obtained the following:


\begin{equation*}
  C_\lambda \times C_\lambda=  [0, 1]\Leftrightarrow   \frac{3-\sqrt{5}}{2} \le\lambda<\frac{1}{2}.
\end{equation*}
In this context, $f_1(C_{1/3}\times C_{1/3})$ is a proper subset of $[0,1]$. Athreya, Reznick and Tyson \cite{ART} showed that $f_2(C_{1/3}\times C_{1/3})=[0,1]$. However, when $k$ in $f_k(x,y)=x^ky$ is large the situation becomes more involved. In the following we provide a characterization of the case $k\in\nn_{\ge2}=\set{n\in\nn:n\ge2}$.

\begin{theorem}\label{theorem-2} Given $k\in \nn_{\ge2}$. If $\lambda_k\le \lambda<\frac{1}{2}$, then $\set{x^ky:x,y\in C_\lambda}$ is exactly interval $[0,1]$, where $\lambda_k$ is the appropriate root of  $(k-1)\lambda^k+(2k+2)\lambda-k-1=0$. Moreover, the sequence $(\lambda_k)$ is increasing to $1/2$ as $k\to \infty$.

\end{theorem}

The rest of our paper is organized as follows. In  section  \ref{S2} we prove Theorem \ref{theorem-1}, and we give the proof of Theorem \ref{theorem-2} in section \ref{S3}. We end the last
section with some remarks.

\section{Proof of Theorem \ref{theorem-1} }\label{S2}

We first introduce some notations. We call $\Gamma_{\sigma}([0,1])=(\Gamma_{i_1}\circ\cdots \circ \Gamma_{i_n})([0,1])$ a basic interval of rank $n$, where $\sigma=i_1\cdots i_n\in\set{0,1}^n$.
Define
\begin{equation*}
\mathcal I_n=\set{I: \text{$I$ is a basic interval of rank $n$}}.
\end{equation*}
For any $I\in \mathcal I_n$, let $I'$ be the union of all basic intervals of rank $n+1$ contained in $I$,
 and denote by $\mathcal I_k^I$ the set of basic intervals of rank $k\ge n$ contained in $I$. 
Now we give a technical lemma: 
\begin{lemma}\label{lemma1}
 Suppose that $\phi: \mathbb{R}^2 \rightarrow \mathbb{R}$ is a continuous function and $\set{F_n}$ be a decreasing sequence of nonempty compact subsets  in $\RR^2$,  then
\begin{equation*}
 \phi\left(\bigcap_{n=1}^\infty F_n \right)=\bigcap_{n=1}^\infty \phi(F_n).
\end{equation*}

\end{lemma}

\begin{proof}
We establish the equality by proving
\begin{equation*}
 \phi\left( \bigcap_{n=1}^\infty F_n \right) \subseteq \bigcap_{n=1}^\infty  \phi(F_n) \qtq{and} \bigcap_{n=1}^\infty  \phi(F_n) \subseteq  \phi\left( \bigcap_{n=1}^\infty F_n \right).
\end{equation*}
Let $ y \in \phi\left( \bigcap_{n=1}^\infty F_n \right) $. Then there exists $ x \in \bigcap_{n=1}^\infty F_n $ such that $ \phi(x) = y $. Since $ x \in F_n $ for all $ n \in \mathbb{N} $, it follows that $ y = \phi(x) \in \phi(F_n) $ for all $ n \in \mathbb{N} $. Therefore, $ y \in \bigcap_{n=1}^\infty \phi(F_n) $.

On the other hand, suppose that $ y \in \bigcap_{n=1}^\infty \phi(F_n) $. For each \( n \in \mathbb{N} \), define
\begin{equation*}
  E_n = F_n \cap \phi^{-1}(\{y\}).
\end{equation*}
Since $ y \in \phi(F_n)$, the set $ E_n $ is nonempty. The continuity of $ \phi $ implies that $ \phi^{-1}(\{y\}) $ is closed, and since $ F_n $ is compact, it follows that $ E_n $ is compact. Furthermore, the sequence $\{E_n\}$ is decreasing because \(\{F_n\}\) is decreasing:
\begin{equation*}
E_{n+1} = F_{n+1} \cap \phi^{-1}(\{y\}) \subseteq F_n \cap \phi^{-1}(\{y\}) = E_n.
\end{equation*}
By the Cantor's intersection theorem, we have $\bigcap_{n=1}^\infty E_n \neq \varnothing$.
Take any \( x \in \bigcap_{n=1}^\infty E_n \). Then
\begin{itemize}
    \item \( x \in F_n \) for all \( n \in \mathbb{N} \), hence \( x \in \bigcap_{n=1}^\infty F_n \);
    \item \( x \in \phi^{-1}(\{y\}) \), hence $ \phi(x) = y $.
\end{itemize}
This shows that $ y \in \phi\left( \bigcap_{n=1}^\infty F_n \right) $.

\end{proof}






In this section we always assume that $f(x, y) = xy$. To prove Theorem \ref{theorem-1} we need several lemmas.

\begin{lemma}\label{lemma2}
Let $ I = [a, a + \lambda^n] $, $ J = [b, b + \lambda^n] $ belong to $ \mathcal{I}_n $ for some $ n \ge 1 $, with $a \le b $. Suppose
\begin{equation*}
\frac{3 - \sqrt{5}}{2} \le \lambda < \frac{1}{2}
\quad \text{and} \quad
b \le \frac{a + \lambda^{n+1}}{1 - 2\lambda}.
\end{equation*}
Then
\begin{equation*}
 f(I, J) = f(I', J').
\end{equation*}

\end{lemma}

\begin{proof}
The sets $ I' $ and $ J' $ are given by
\begin{align*}
I' &= [a, a + \lambda^{n+1}] \cup [a + \lambda^n - \lambda^{n+1}, a + \lambda^n], \\
J' &= [b, b + \lambda^{n+1}] \cup [b + \lambda^n - \lambda^{n+1}, b + \lambda^n].
\end{align*}
The image $f(I', J') $ decomposes into the following four intervals:
\begin{align*}
I_1 &= \left[ ab,\; (a + \lambda^{n+1})(b + \lambda^{n+1}) \right], \\
I_2 &= \left[ ab + a(1 - \lambda)\lambda^n,\; (a + \lambda^{n+1})(b + \lambda^n) \right], \\
I_3 &= \left[ ab + b(1 - \lambda)\lambda^n,\; (a + \lambda^n)(b + \lambda^{n+1}) \right], \\
I_4 &= \left[ (a + \lambda^n - \lambda^{n+1})(b + \lambda^n - \lambda^{n+1}),\; (a + \lambda^n)(b + \lambda^n) \right].
\end{align*}
Observe that $ f(I, J) = \left[ ab,\; (a + \lambda^n)(b + \lambda^n) \right] $. To establish the equality $f(I, J) = f(I', J') $, it is enough to verify that the intervals $ I_1, \dots, I_4 $ cover $ [ab, (a + \lambda^n)(b + \lambda^n)] $. This is ensured by the three inequalities:
\begin{equation*}
  r_1 \ge l_2, \quad r_2 \ge l_3, \quad r_3 \ge l_4,
\end{equation*}
where $ I_i = [l_i, r_i] \) for \( i = 1, \dots, 4 $.

\noindent \text{(i) Verification of $ r_1 \ge l_2 $.}
We require:
\begin{equation*}
 (a + \lambda^{n+1})(b + \lambda^{n+1}) \ge ab + a(1 - \lambda)\lambda^n.
\end{equation*}
This is equivalent to $\lambda^{n+2} + (2a + b)\lambda - a \ge 0$.
Since $ a \le b $ and $\lambda > \frac{1}{3}$   we have
\begin{equation*}
 \lambda^{n+2} + (2a + b)\lambda - a > (2a + b)\lambda - a\ge a(3\lambda-1)>0.
\end{equation*}


\noindent \text{(ii) Verification of \( r_3 \ge l_4 \).}
We require:
\begin{equation*}
 (a + \lambda^n)(b + \lambda^{n+1}) \ge (a + \lambda^n - \lambda^{n+1})(b + \lambda^n - \lambda^{n+1}).
\end{equation*}
This is equivalent to:
\begin{equation*}
 (2a + b)\lambda + 3\lambda^{n+1} -a - \lambda^n - \lambda^{n+2} \ge 0.
\end{equation*}
Note that
\begin{equation*}
 - \lambda^n + 3\lambda^{n+1} - \lambda^{n+2} = \lambda^n(3\lambda - \lambda^2 - 1).
\end{equation*}
The assumption $ \lambda \ge \frac{3 - \sqrt{5}}{2} $ implies $3\lambda - \lambda^2 - 1 \ge 0 $, together with $(2a + b)\lambda - a\ge0$ we have $r_3\ge l_4$.

\noindent \text{(iii) Verification of \( r_2 \ge l_3 \).}
We require:
\begin{equation*}
  (a + \lambda^{n+1})(b + \lambda^n) \ge ab + b(1 - \lambda)\lambda^n.
\end{equation*}
This can be checked by the assumed inequality $b \le \frac{a + \lambda^{n+1}}{1 - 2\lambda}$.

 Since all three inequalities are satisfied, the intervals \( I_1, \dots, I_4 \) cover \( [ab, (a + \lambda^n)(b + \lambda^n)] \) without gaps, and the lemma follows.
\end{proof}

 Given two basic intervals $I,J\in\mathcal I_n$, Jiang, Kong, Li and Wang \cite{JKLW} introduced the \emph {double covering set} $D(I,J)$, defined as the set of points $x \in f(I,J)$ for which there exist two distinct pairs $(I_p, J_p), (I_q, J_q) \in \mathcal{I}_k^I \times \mathcal{I}_k^J$ (with $k > n$) such that $x$ lies in the interior of both $f(I_p, J_p)$ and $f(I_q, J_q)$.
The following result shows that the interior of $f(I,J)$ is double covering set under some conditions.

\begin{lemma}\label{lemma3}
Let $ I = [a, a + \lambda^n] $, $ J = [b, b + \lambda^n] \in \mathcal{I}_n $ for some $ n \ge 1 $ with $ a \le b $. If
\begin{equation}\label{condition-1}
 \frac{3 - \sqrt{5}}{2}<\lambda < \frac{1}{2}
\quad \text{and} \quad
b \le \frac{a \lambda}{1 - 2\lambda},
\end{equation}
then
\begin{equation*}
(L_n(a,b),\; R_n(a,b)) \subset D(I,J),
\end{equation*}
where
\begin{equation*}
 L_n(a,b) = ab + a(1 - \lambda)\lambda^n, \quad
R_n(a,b) = (a + \lambda^n)(b + \lambda^{n+1}).
\end{equation*}
\end{lemma}

\begin{proof}
We adopt the notation from the proof of Lemma~\ref{lemma2}:
\begin{align*}
[l_1, r_1] &= \left[ ab,\; (a + \lambda^{n+1})(b + \lambda^{n+1}) \right], \\
[l_2, r_2] &= \left[ ab + a(1 - \lambda)\lambda^n,\; (a + \lambda^{n+1})(b + \lambda^n) \right], \\
[l_3, r_3] &= \left[ ab + b(1 - \lambda)\lambda^n,\; (a + \lambda^n)(b + \lambda^{n+1}) \right], \\
[l_4, r_4] &= \left[ (a + \lambda^n - \lambda^{n+1})(b + \lambda^n - \lambda^{n+1}),\; (a + \lambda^n)(b + \lambda^n) \right].
\end{align*}
We show that $ l_3< r_1 $ and $ l_4 < r_2 $. For  $ l_3< r_1 $,
we require:
\begin{equation*}
ab + b(1 - \lambda)\lambda^n < (a + \lambda^{n+1})(b + \lambda^{n+1}).
\end{equation*}
This is implied by the assumption \eqref{condition-1}.
For the second statement $ l_4 < r_2 $ we require:
\begin{equation*}
(a + \lambda^n - \lambda^{n+1})(b + \lambda^n - \lambda^{n+1}) < (a + \lambda^{n+1})(b + \lambda^n).
\end{equation*}
This is equivalent to:
\begin{equation*}
b(1 - 2\lambda) <a\lambda - (1 - 3\lambda + \lambda^2)\lambda^n.
\end{equation*}
Since $ 1 - 3\lambda + \lambda^2 < 0 $ for $ \lambda >\frac{3 - \sqrt{5}}{2} $, the right-hand side is at least $ a\lambda $, so  $ l_4 < r_2 $ follows from \eqref{condition-1}. 
 We conclude that the interval $(l_2, r_3) = (L_n(a,b),\; R_n(a,b))$ is contained in $ D(I,J)$, as required.

\end{proof}

Combined with the fact that $ l_2 \le l_3,r_2 \le r_3 $ (which follows from $a \le b $), Lemma \ref{lemma3} implies that the overlapping pattern of $[l_i,r_i]$ coincides with the description in Figure 1.

\begin{figure}[h]\label{figure1}
\centering
\begin{tikzpicture}[scale=5]

\draw(-1.2,.4)node[above]{\scriptsize $l_1$}--(-.4,.4)node[below]{\scriptsize $r_1$};

\draw(-.8,.48)node[above]{\scriptsize $l_2$}--(0,.48)node[above]{\scriptsize $r_2$};

\draw(-.6,.32)node[below]{\scriptsize $l_3$}--(.2,.32)node[below]{\scriptsize $r_3$};

\draw(-.2,.4)node[left]{\scriptsize $l_4$}--(.6,.4)node[below]{\scriptsize $r_4$};

\end{tikzpicture}\caption{The overlapping pattern of $[l_i,r_i],i=1,2,3,4$.}
\end{figure}
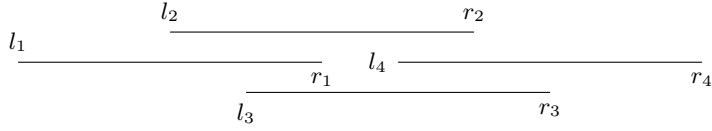

\begin{lemma}\label{pro-1}
 Let $ I = [a, a + \lambda^n]  , J = [b, b + \lambda^n] \in \mathcal I_n$ for some $n\ge1$ with $a\le b$. If $\frac{3-\sqrt{5}}{2}< \lambda<\frac{1}{2}$ and
 \begin{equation}\label{p1}
  \frac{(1-\lambda-\lambda^2)(a+\lambda^n)}{\lambda}\le b<b+\lambda^n \le  \frac{a\lambda}{1-2\lambda}.
 \end{equation}
Then  we have $D(I,J)=\mathrm{int}( f(I,J))$.
\end{lemma}

\begin{proof}
We divide the proof into two steps.

Step (I). We prove that
\begin{equation*}
 (l_1,r_3)=(ab,(a+\lambda^{n})(b+\lambda^{n+1}))\subset D(I,J).
\end{equation*}
Let $I_k = [a, a + \lambda^k]$ and $J_k = [b, b + \lambda^k]$ for $k \geq n$.
Since $I_k \subset I$ and $J_k \subset J$, it follows that
\begin{equation*}
\bigcup_{k=n}^{\infty} D(I_k, J_k) \subset D(I, J).
\end{equation*}
By \eqref{p1} for each $k \geq n$, we have
\begin{equation*}
b < b + \lambda^k \leq b + \lambda^n \leq \frac{a\lambda}{1 - 2\lambda}.
\end{equation*}
This implies $(L_k(a,b), R_k(a,b)) \subset D(I_k, J_k)$ by Lemma \ref{lemma3}.
Consequently,
\begin{equation*}
\bigcup_{k=n}^{\infty} (L_k(a,b), R_k(a,b)) \subset D(I, J).
\end{equation*}

On the other hand, by \eqref{p1} it follows that
\begin{align*}
  L_k(a,b)-R_{k+1}(a,b)&=ab+a(1- \lambda)\lambda^{k}-(a+\lambda^{k+1})(b+\lambda^{k+2})\\
                       &=-\lambda^{k}[\lambda^{k+3}+a(\lambda+\lambda^2-1)+b\lambda]< 0.
\end{align*}
Furthermore, $L_k(a,b)\rightarrow ab$ and $R_{k}(a,b)\rightarrow ab$ as $k\to \infty$. So we have \begin{equation*}
\bigcup_{k=n}^{\infty}(L_k(a,b),R_k(a,b))=(ab,(a+\lambda^{n})(b+\lambda^{n+1}))\subset D(I,J).
\end{equation*}

Step (II). Let $k\ge n$ and
\begin{align*}
  &\tilde{I_k}=[a+\lambda^n-\lambda^k,a+\lambda^n]=[a_k,a_k+\lambda^k],\\
  &\tilde{J_k}=[b+\lambda^n-\lambda^k,b+\lambda^n]=[b_k,b_k+\lambda^k].
\end{align*}
Then $\tilde{I_k}\subset I,\tilde{J_k}\subset J$ and  $\cup_{k=n}^{\infty}D(\tilde{I_k},\tilde{J_k})\subset D(I,J)$. It follows from \eqref{p1} that
\begin{equation*}
\frac{(1-\lambda-\lambda^2)(a_k+\lambda^k)}{\lambda}=\frac{(1-\lambda-\lambda^2)(a+\lambda^n)}{\lambda}\le b\le b_k
\end{equation*}
and
\begin{equation*}
 b_k+\lambda^k=b+\lambda^n \le  \frac{a\lambda}{1-2\lambda}\le \frac{a_k\lambda}{1-2\lambda}.
\end{equation*}
So using the result in Part (I) for $\tilde{I_k},\tilde{J_k}$ we get
\begin{equation*}
\bigcup_{k=n}^{\infty}(a_kb_k,(a_k+\lambda^{k})(b_k+\lambda^{k+1}))\subset\bigcup_{k=n}^{\infty}D(\tilde{I_k},\tilde{J_k})\subset D(I,J).
\end{equation*}
Finally, one can show that $(a_k+\lambda^{k})(b_k+\lambda^{k+1})> a_{k+1}b_{k+1}$. This can be verified by
the following calculation:
\begin{align*}
  &(a_k+\lambda^{k})(b_k+\lambda^{k+1})-a_{k+1}b_{k+1}\\
                       =&(a+\lambda^{n})(b+\lambda^{n}+\lambda^{k+1}-\lambda^{k})-(a+\lambda^{n}-\lambda^{k+1})(b+\lambda^{n}-\lambda^{k+1})\\
         =&\lambda^{k}[(2a+b)\lambda-a+3\lambda^{n+1}-\lambda^n-\lambda^{k+2}].
\end{align*}
Using \eqref{p1} again we have $(2a+b)\lambda-a\ge0$ and $3\lambda-1-\lambda^2>0$. Since $\lim\limits_{k\to\infty}(a_k+\lambda^{k})(b_k+\lambda^{k+1})=(a+\lambda^{n})(b+\lambda^{n})$, therefore
\begin{equation*}
  \bigcup_{k=n}^{\infty}(a_kb_k,(a_k+\lambda^{k})(b_k+\lambda^{k+1}))=(ab,(a+\lambda^{n})(b+\lambda^{n}))\subset D(I,J),
\end{equation*}
and we conclude that $D(I,J)=\mathrm{int}( f(I,J))$.
\end{proof}

We recall that $S_t=\set{(x,y)\in C_\lambda\times C_\lambda:xy=t}$ and the first main result:

\begin{theorem}\label{theorem-11} If $0.4302\le \lambda<\frac{1}{2}$, then $S_t$ has the cardinality of the continuum for any $0<t<1$.
\end{theorem}

\begin{proof}
The proof is divided into two steps.

Step (I). Let $ I = [a, a + \lambda^n]  , J = [b, b + \lambda^n] \in \mathcal I_n$ for some $n\ge1$ with $a\le  b$. We show that $S_t$ is of cardinality
continuum for any $ab<t<(a+\lambda^{n})(b+\lambda^{n})$.

Given $ab<t<(a+\lambda^{n})(b+\lambda^{n})$. Note that for any $I_k=[a_k,a_k+\lambda^{k}],J_k=[b_k,b_k+\lambda^{k}],k\ge n$, where $a\le a_k \le a+\lambda^n-\lambda^k$ and $b\le b_k \le b+\lambda^n-\lambda^k$, if \eqref{p1} holds for $I$ and $J$ then we have
\begin{equation*}
\frac{(1-\lambda-\lambda^2)(a_k+\lambda^k)}{\lambda}\le \frac{(1-\lambda-\lambda^2)(a+\lambda^n)}{\lambda}\le b\le b_k
\end{equation*}
and
\begin{equation*}
 b_k+\lambda^k\le b+\lambda^n \le  \frac{a\lambda}{1-2\lambda}\le \frac{a_k\lambda}{1-2\lambda}.
\end{equation*}
So it follows from Lemma \ref{pro-1} that
\begin{equation}\label{int}
\mathrm{int} (f(I_k,J_k))= D(I_k,J_k).
\end{equation}
We distinguish the following
two cases: (A) $a<b$; (B) $a=b$.

Case (A). $a<b$.  We inductively define a binary tree of basic interval pairs $(I_{i_1\ldots i_n},J_{i_1\ldots i_n})$ where $i_n\in \set{p,q},n\in\nn$.  We may depict this bifurcation as shown in Figure 2.

By Lemma \ref{pro-1} there exist two different interval pairs $(I_{p},J_{p}),(I_{q},J_{q})\in \mathcal I_{k}^I\times \mathcal I_{k}^J$ such that $t\in \mathrm{int} (f(I_{p},J_{p}))$ and $t\in \mathrm{int} (f(I_{q},J_{q}))$, where $(I_{p},J_{p})$ and $(I_{q},J_{q})$
are the first and second interval pairs selected from $\mathcal I_{k}^I\times \mathcal I_{k}^J$, respectively.
 Suppose that $i_1,i_2,\ldots,i_{n}$ have been defined, since $t\in D(I_{i_1 i_2 \ldots i_{n}},J_{i_1 i_2 \ldots i_{n}})$ by \eqref{int}, there exist two different basic intervals $I_{i_1 i_2\ldots i_{n}p},I_{i_1 i_2 \ldots i_{n}q}\in \mathcal I^{I_{i_1i_2\ldots i_{n}}}$, and $J_{i_1i_2\ldots i_{n}p},J_{i_1i_2\ldots i_{n}q}\in \mathcal I^{J_{i_1i_2\ldots i_{n}}}$ such that $t\in \mathrm{int} ( f(I_{i_1i_2\ldots i_{n}p},J_{i_1i_2\ldots i_{n}p}))$ and $t\in \mathrm{int} (f(I_{i_1i_2\ldots i_{n}q},J_{i_1i_2\ldots i_{n}q}))$. Repeat this process again and continue. Since $(I_{i_1i_2\ldots i_{n}}),(J_{i_1i_2\ldots i_{n}})$ are decreasing sequences of nonempty compact set, for any sequence $\mathbf{c}=(i_{n})_{n=1}^\infty \in\set{p,q}^\nn$ let
\begin{equation*}
x_\mathbf{c}=\bigcap_{n=1}^\infty I_{i_1i_2\ldots i_{n}}\in C_\lambda,y_\mathbf{c}=\bigcap_{n=1}^\infty J_{i_1i_2\ldots i_{n}}\in C_\lambda.
\end{equation*}
Observe that  $t\in f(I_{i_1i_2\ldots i_{n}},J_{i_1i_2\ldots i_{n}})$ for any $n\in\nn$, so it follows from Lemma \ref{lemma1} that
\begin{equation*}
  t\in \bigcap_{n=1}^\infty f(I_{i_1i_2\ldots i_{n}},J_{i_1i_2\ldots i_{n}})= f\left(\bigcap_{n=1}^\infty I_{i_1i_2\ldots i_{n}},\bigcap_{n=1}^\infty J_{i_1i_2\ldots i_{n}}\right).
\end{equation*}
Then we have $x_\mathbf{c}y_\mathbf{c}=t$, which implies that
\begin{equation*}
 \set{(x_\mathbf{c},y_\mathbf{c}):x_\mathbf{c}y_\mathbf{c}=t,\mathbf{c}\in\set{p,q}^\nn}\subset S_t.
\end{equation*}
So $S_t$ is of cardinality
continuum.

\begin{figure}[h]\label{figure2}
\centering
\begin{tikzpicture}[scale=5]

\draw(-1.2,.4)node[above]{\scriptsize $(I,J)$}--(-.8,.4);

\draw(-.8,.6)node[above]{\scriptsize $(I_{p},J_{p})$}--(-.8,.2)node[below]{\scriptsize $(I_{q},J_{q})$};

\draw(-.8,.6)--(-.5,.6);

\draw(-.5,.7)node[above]{\scriptsize $(I_{pp},J_{pp})$}--(-.5,.5)node[below]{\scriptsize $(I_{pq},J_{pq})$};

\draw(-.5,.7)--(-.3,.7)node[right]{$\cdots$};

\draw(-.5,.5)--(-.3,.5)node[right]{$\cdots$};

\draw(-.8,.2)--(-.2,.2);

\draw(-.2,.3)node[above]{\scriptsize $(I_{qp},J_{qp})$}--(-.2,.1)node[below]{\scriptsize $(I_{qq},J_{qq})$};

\draw(-.2,.3)--(0,.3)node[right]{$\cdots$};

\draw(-.2,.1)--(0,.1)node[right]{$\cdots$};

\end{tikzpicture}\caption{Binary tree of basic interval pairs.}
\end{figure}
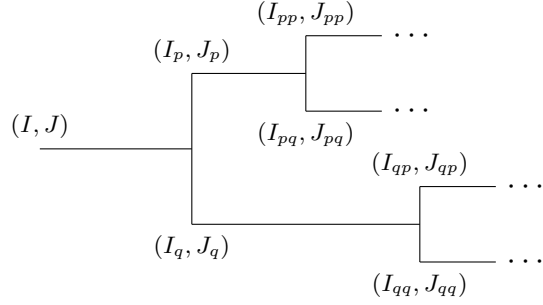

Case (B). $a=b$. It follows from Lemma \ref{pro-1} that there exist two different interval pairs $(I_{p},J_{p}),(I_{q},J_{q})\in \mathcal I_{k}^I\times \mathcal I_{k}^J$ such that $t\in \mathrm{int} (f(I_{p},J_{p}))$ and $t\in \mathrm{int} (f(I_{q},J_{q}))$. We distinguish two cases: (i) $I_{p}\neq J_{p}$ or $I_{q}\neq J_{q}$; (ii) $I_{p}=J_{p}$ and $I_{q}=J_{q}$.

(i) $I_{p}\neq J_{p}$ or $I_{q}\neq J_{q}$. We assume that $I_{p}\neq J_{p}$. By Lemma \ref{pro-1} there exist two different interval pairs $(I_{pp},J_{pp}),(I_{pq},J_{pq})\in \mathcal I_{k}^I\times \mathcal I_{k}^J$ such that $t\in \mathrm{int} (f(I_{pp},J_{pp}))$ and $t\in \mathrm{int} (f(I_{pq},J_{pq}))$. 
In view the proof of Case (A), we can inductively define a binary tree of basic interval pairs $(I_{i_1\ldots i_n},J_{i_1\ldots i_n})$ where $i_n\in \set{p,q},n\in\nn$. So $S_t$ is of cardinality continuum. For the case $I_{q}\neq J_{q}$ the proof is similar.

(ii) $I_{p}=J_{p}$ and $I_{q}=J_{q}$. Again by Lemma \ref{pro-1}, we obtain two different interval pairs $(I_{pp},J_{pp}),(I_{pq},J_{pq})\in \mathcal I_{k}^I\times \mathcal I_{k}^J$ such that $t\in \mathrm{int} (f(I_{pp},J_{pp}))$ and $t\in \mathrm{int} (f(I_{pq},J_{pq}))$, and two different interval pairs $(I_{qp},J_{qp}),(I_{qq},J_{qq})\in \mathcal I_{k}^I\times \mathcal I_{k}^J$ such that $t\in \mathrm{int} (f(I_{qp},J_{qp}))$ and $t\in \mathrm{int} (f(I_{qq},J_{qq}))$.  If any of the four cases $I_{pp}\neq J_{pp}, I_{pq}\neq J_{pq}, I_{qp}\neq J_{qp}, I_{qq}\neq J_{qq}$ occurs then the proof is complete.  Note that $I_{i_1 i_2 \ldots i_{n}}=J_{i_1 i_2 \ldots i_{n}}$ cannot happen uncountably many times, because when $x=\cap_{n=1}^\infty I_{i_1 i_2\ldots i_{n}}=\cap_{n=1}^\infty J_{i_1 i_2 \ldots i_{n}}=y$, $xy=t$ has the unique solution. This implies that $I_{i_1 i_2 \ldots i_{n}}=J_{i_1 i_2 \ldots i_{n}}$ occurs in every branch is impossible. So there exists two intervals $I_{i_1 i_2 \ldots i_{n}}\neq J_{i_1 i_2 \ldots i_{n}}$ in some branch, by the argument in Case (A) we complete the proof.

Step (II). Now we give some basic intervals and verify that they  satisfy \eqref{p1}. First we choose $ I_{3,1} = J_{3,1}=[1-\lambda^3, 1] $, then we obtain 
\begin{align*}
   \frac{1-\lambda-\lambda^2}{\lambda}\le 1-\lambda^3 (\lambda\ge 0.4258), 1 \le  \frac{(1-\lambda^3)\lambda}{1-2\lambda}   (\lambda\ge 0.3377).
\end{align*}
Hence,  \eqref{p1} is satisfied.  Let $ I_{3,2} =[1-\lambda^2, 1-\lambda^2+\lambda^3],J_{3,2}=[1-\lambda^3, 1]$, then we have
\begin{align*}
   \frac{(1-\lambda-\lambda^2)(1-\lambda^2+\lambda^3)}{\lambda}\le 1-\lambda^3 (\lambda\ge0.4094), 1 \le  \frac{(1-\lambda^2)\lambda}{1-2\lambda}  (\lambda\ge0.3474).
\end{align*}
Therefore,  \eqref{p1} is again satisfied. We show the rest intervals and the range of $\lambda$ satisfying \eqref{p1} in Table 1.

\begin{table}[!ht]
    \begin{center}
    \label{tab:first_table}
    \begin{tabular}{|c|c|c|}
        \hline
        \textbf{$I_{3,n},n=3,4,5,6$} & \textbf{$J_{3,n},n=3,4,5,6$} & \text{The range of $\lambda$ } \\
        \hline
        $[1-\lambda^2, 1-\lambda^2+\lambda^3]$ & $[1-\lambda^2, 1-\lambda^2+\lambda^3]$ & $\lambda\ge0.4274$ \\
        \hline
       $[1-\lambda+\lambda^2-\lambda^3, 1-\lambda+\lambda^2]$ & $ [1-\lambda^2, 1-\lambda^2+\lambda^3]$ & $\lambda\ge0.3993$ \\
        \hline
        $[1-\lambda+\lambda^2-\lambda^3, 1-\lambda+\lambda^2]$ & $[1-\lambda+\lambda^2-\lambda^3, 1-\lambda+\lambda^2]$ & $\lambda\ge0.4302$  \\
        \hline
        $[1-\lambda, 1-\lambda+\lambda^3]$ & $[1-\lambda+\lambda^2-\lambda^3, 1-\lambda+\lambda^2]$ & $\lambda\ge0.4076$  \\
        \hline
    \end{tabular}

    \end{center}
     \caption{The selected basic intervals and the corresponding ranges of $\lambda$ .}
\end{table}

Furthermore, we have
\begin{equation*}
  \bigcup_{n=1}^6\mathrm{int}(f(I_{3,n},J_{3,n}))=((1 -\lambda)(1-\lambda+\lambda^2-\lambda^3),1).
\end{equation*}
This can
be verified by the following calculation:
\begin{align*}
 &1-\lambda^2+\lambda^3> (1-\lambda^3)^2(\lambda\ge0.3377),\\
 &(1-\lambda^2+\lambda^3)^2>(1-\lambda^2)(1-\lambda^3) (\lambda\ge 0.3290),\\
 &(1-\lambda+\lambda^2)(1-\lambda^2+\lambda^3)>(1-\lambda^2)^2 (\lambda\ge0.4198),\\
 &(1-\lambda+\lambda^2)^2>(1-\lambda+\lambda^2-\lambda^3)(1-\lambda^2)(\lambda\ge0.4084),\\
&(1-\lambda+\lambda^3)(1-\lambda+\lambda^2)>(1-\lambda+\lambda^2-\lambda^3)^2 (\lambda\ge0.3391).
\end{align*}
Based on the above reasoning we obtain $S_t$ is of cardinality
continuum for any $t\in ((1-\lambda)(1-\lambda+\lambda^2-\lambda^3),1)$. Since for any $z\in C_\lambda$ we have $\lambda z\in C_\lambda$, so if $x,y\in C_\lambda$ satisfying $xy=t$ for some $t\in ((1-\lambda)(1-\lambda+\lambda^2-\lambda^3),1)$, then $\lambda x,y\in C_\lambda$ satisfy $\lambda xy=\lambda t$ for $\lambda t\in (\lambda(1-\lambda)(1-\lambda+\lambda^2-\lambda^3),\lambda)$. Hence the  conclusion  is deduced from the fact that $\lambda>(1-\lambda)(1-\lambda+\lambda^2-\lambda^3)$ with $\lambda>0.4084$.

\end{proof}
\begin{remark}
 In Case (B)(i) of the proof of Theorem \ref{theorem-11}, if the case $I_{p}$  is located to the right of $J_{p}$ occurs in some branch, the symmetry of $f(x,y)$ yields $f(I_{p},J_{p})=f(J_{p},I_{p})$. We choose $(J_{p},I_{p})$ and redefine it as $(I_{p},J_{p})$, then $(I_{p},J_{p})$ satisfies the condition of Lemma \ref{pro-1}. The case $I_{q}$  is located to the right of $J_{q}$ is analogous.
\end{remark}

\section{Proof of Theorem \ref{theorem-2} }\label{S3}
Given $k\in \nn_{\ge2}$ and $ f_k(x, y) = x^ky$, denote by $\lambda_k$  be the appropriate root of
\begin{equation*}
 (k-1)\lambda^k+2(k+1)\lambda-k-1=0.
\end{equation*}
Recall  that $\mathcal I_n$ is the set of all basic intervals of rank $n$. For any $I\in \mathcal I_n$,  $I'$ is the union of all basic intervals of rank $n+1$ contained in $I$. This section we focus on the case that the element in $\mathcal I_n$ is the subset of $[1-\lambda,1]$.

\begin{lemma}\label{lemma31} Let  $ I = [a, a + \lambda^n]  , J = [b, b + \lambda^n] \in \mathcal I_n$ for some $n\ge1$. If $\lambda_k\le \lambda<\frac{1}{2}$ and
\begin{equation}\label{e31}
 \frac{(1-2\lambda)(a+k\lambda^n)}{k\lambda} \le b \le  \frac{a}{k(1-2\lambda)}.
\end{equation}
Then $ f_k(I, J)= f_k(I', J')$.
\end{lemma}
\begin{proof}
We observe that $ f_k(I, J)=[a^kb,(a+\lambda^{n})^k(b+\lambda^{n})]$ and $f_k(I', J')$ can be divided into 4 sub-intervals
\begin{align*}
&[u_1,v_1] = \left[ a^kb, (a+\lambda^{n+1})^k(b+\lambda^{n+1})   \right],\\
&[u_2,v_2]  = \left[ a^kb+a^k\lambda^{n}(1- \lambda), (a+\lambda^{n+1})^k(b+\lambda^{n})  \right],\\
&[u_3,v_3]  = \left[ (a+\lambda^{n}- \lambda^{n+1})^kb, (a+\lambda^{n})^k(b+\lambda^{n+1}) \right],\\
&[u_4,v_4]  = \left[ (a+\lambda^{n}- \lambda^{n+1})^k(b+\lambda^{n}- \lambda^{n+1}), (a+\lambda^{n})^k(b+\lambda^{n}) \right].
\end{align*}
It suffices to establish the inequalities $v_1\ge u_2,v_2\ge u_3$ and $v_3\ge u_4$ under the condition \eqref{e31}.

\noindent \text{(i) Verification of $ v_1\ge u_2 $.}
First, $b\ge \frac{a(1-2\lambda)}{k\lambda}$ implies
\begin{align*}
 (a^k+ka^{k-1}\lambda^{n+1})(b+\lambda^{n+1})\ge a^k[b+(1- \lambda)\lambda^{n}].
\end{align*}
It follows that
\begin{equation*}
  (a+\lambda^{n+1})^k(b+\lambda^{n+1})\ge a^k[b+(1- \lambda)\lambda^{n}].
\end{equation*}
\noindent \text{(ii) Verification of $ v_3\ge u_4 $.} Write $c=(a+\lambda^{n})^k-a^k$. Since $ b\ge  \frac{(1-2\lambda)(a+k\lambda^n)}{k\lambda}$ and $c>ka^{k-1}\lambda^{n}$, a simple calculation yields
\begin{equation*}
   b\ge \frac{(\frac{a^k\lambda^{n}}{c}+\lambda^n)(1- 2\lambda)}{\lambda}=\frac{(a^k\lambda^{n}+c\lambda^n)(1-2\lambda)}{c\lambda}.
\end{equation*}
It follows that
\begin{align*}
  (a+\lambda^{n})^k(b+\lambda^{n+1}) &\ge [a^k+c(1-\lambda)][b+\lambda^{n}(1- \lambda)]\\
                                    &\ge [a+\lambda^{n}(1- \lambda)]^k[b+\lambda^{n}(1- \lambda)].
\end{align*}
\noindent \text{(iii) Verification of $ v_2\ge u_3 $.} To establish $v_2\ge u_3$, let $c_{i,k}=C_k^ia^{k-i}\lambda^{ni}(1\le i\le k)$, where $C_k^i=k(k-1)\cdots (k-i+1)/i!$. We require
\begin{align*}
   (a+\lambda^{n+1})^k(b+\lambda^{n})  \ge [a+(1- \lambda)\lambda^{n}]^kb,
\end{align*}
that is
\begin{equation*}
(a+\lambda^{n+1})^k\lambda^{n}+\sum_{i=1}^kc_{i,k}\lambda^{i}b \ge \sum_{i=1}^kc_{i,k}(1-\lambda)^ib.
\end{equation*}
We point that $(1-\lambda)^{i+1}-\lambda^{i+1}\le (1-\lambda)^i-\lambda^i\le 1-2\lambda$ for any $i\in\zz^+$.
To prove  $(a+\lambda^{n+1})^k\lambda^{n}+\sum_{i=1}^kc_{i,k}\lambda^{i}b \ge \sum_{i=1}^kc_{i,k}(1-\lambda)^ib$, it suffices to show that
\begin{equation*}
 \left(a^k+\sum_{i=1}^kc_{i,k}\lambda^i\right)\lambda^{n} \ge \sum_{i=1}^kc_{i,k}(1-2\lambda)b.
\end{equation*}
By the assumption $bk(1-2\lambda)\le a$ we obtain $a^k \lambda^{n}\ge c_{1,k}(1-2\lambda)b$. It remains to prove that $\sum_{i=2}^kc_{i,k}(1-2\lambda)b\le \sum_{i=1}^{k}c_{i,k}\lambda^{n+i}$. In fact we have
\begin{equation}\label{e32}
\sum_{i=2}^kc_{i,k}(1-2\lambda)b<  \sum_{i=2}^kc_{i,k}(1-2\lambda) < \sum_{i=1}^{k-1}c_{i,k}\lambda^{n+i}.
\end{equation}
This will be proved by induction on $k$. When $k=2$, since $a\ge 1-\lambda> \lambda$ and $\lambda\ge \lambda_2>(\sqrt{3}-1)/2$, we have $1-2\lambda\le 2\lambda^2< 2a\lambda$. When $k=3$, by assumption $\lambda\ge \lambda_3>0.41$, we have $1<(3-\lambda)\lambda \le (a+2)\lambda$ and $1<3\lambda^3+2\lambda< 3a\lambda^2+2\lambda$. So $3a<3a^2\lambda+6a\lambda$ and $\lambda^n<3a\lambda^{n+2}+2\lambda^{n+1}$, it follows that $ (1-2\lambda)(3a+\lambda^n)< 3a^2\lambda+3a\lambda^{n+2}$. We rewrite the second inequality in \eqref{e32} and assume that
 \begin{equation}\label{e36}
   (1-2\lambda)\sum_{j=1}^{k-1}C_k^{j+1}a^{k-j-1}\lambda^{nj}< \sum_{j=1}^{k-1}C_k^ja^{k-j}\lambda^{(n+1)j}
 \end{equation}
 holds for some $k\ge3$. For $1\le j\le k-1$ we consider the differences
\begin{align*}
  &C_{k+1}^{j+1}a^{k-j}\lambda^{nj}-C_{k}^{j+1}a^{k-j}\lambda^{nj}, \\
  &C_{k+1}^ja^{k+1-j}\lambda^{(n+1)j}-C_{k}^ja^{k+1-j}\lambda^{(n+1)j}.
\end{align*}
Note that
 \begin{align}\label{e35}
   (1-2\lambda)(C_{k+1}^{j+1}a^{k-j}\lambda^{nj}-C_{k}^{j+1}a^{k-j}\lambda^{nj})<C_{k+1}^ja^{k+1-j}\lambda^{(n+1)j}-C_{k}^ja^{k+1-j}\lambda^{(n+1)j},
 \end{align}
which is equivalent to
\begin{align}\label{e37}
  &k-j+1< ja\lambda^{j}+2(k-j+1)\lambda.
\end{align}
It can be verified that $m\lambda^m>(m+1)\lambda^{m+1}$ for $m\in\zz^+$, thus we obtain $\frac{ja\lambda^{j}}{k+1}+2\lambda \ge\frac{(k-1)a\lambda^{k-1}}{k+1}+2\lambda>\frac{(k-1)\lambda^{k}}{k+1}+2\lambda\ge1$. This implies \eqref{e37}. Take the sum for $j$ from $1$ to $k-1$ on both sides of  \eqref{e35} and combine with \eqref{e36} we have
\begin{equation*}
(1-2\lambda)\sum_{j=1}^{k-1}C_{k+1}^{j+1}a^{k-j}\lambda^{nj}< \sum_{j=1}^{k-1}C_{k+1}^ja^{k+1-j}\lambda^{(n+1)j}.
\end{equation*}
Moreover, it follows from $a(k+1)>\lambda(k+1)>\frac{k-1}{k+1}$ (because $\lambda\ge\lambda_k\ge \lambda_2>0.46$) that $(1-2\lambda)\lambda^{nk}< (k+1)a\lambda^{(n+1)k}$.

Finally, from the above reasoning, it follows that
\begin{align*}
 &(1-2\lambda)\sum_{j=1}^{k}C_{k+1}^{j+1}a^{k-j}\lambda^{nj}=(1-2\lambda)\sum_{j=1}^{k-1}C_{k+1}^{j+1}a^{k-j}\lambda^{nj}+(1-2\lambda)\lambda^{nk}\\
&< \sum_{j=1}^{k-1}C_{k+1}^ja^{k+1-j}\lambda^{(n+1)j}+(k+1)a\lambda^{(n+1)k}=\sum_{j=1}^{k}C_{k+1}^ja^{k+1-j}\lambda^{(n+1)j}.
\end{align*}
This proves \eqref{e36} holds for $ k + 1$, and hence  \eqref{e32}  follows by induction.

\end{proof}

\begin{lemma}\label{lemma32} Let $ I = [a, a + \lambda^n]  , J = [b, b + \lambda^n] \in \mathcal I_n$ for some $n\ge1$. If $\lambda_k\le \lambda<\frac{1}{2}$ and
\begin{equation}\label{e33}
 \frac{(1-2\lambda)(a+k\lambda^n)}{k\lambda} \le b <b+\lambda^n \le  \frac{a}{k(1-2\lambda)}.
\end{equation}
Then $ f_k(I, J)= f_k(C_\lambda \cap I, C_\lambda \cap J)$.
\end{lemma}
\begin{proof}
Given $m\ge n$, let $I_m=[a_0,a_0+\lambda^m]\subset I,J_m=[b_0,b_0+\lambda^m]\subset J$ be two basic intervals. By \eqref{e33} and our assumption $a\le a_0\le a+\lambda^n-\lambda^m$ and $b\le b_0\le b+\lambda^n-\lambda^m$, we obtain
\begin{equation*}
 \frac{(1-2\lambda)(a_0+k\lambda^m)}{k\lambda} \le \frac{(1-2\lambda)(a+k\lambda^n)}{k\lambda}\le b\le b_0,
\end{equation*}
and
\begin{equation*}
 b_0+\lambda^m\le b+\lambda^n\le   \frac{a}{k(1-2\lambda)}\le   \frac{a_0}{k(1-2\lambda)}.
\end{equation*}
It follows from Lemma \ref{lemma31} that $f_k(I_m,J_m)=f_k(\cup \mathcal I_{m+1}^{I_m},\cup \mathcal I_{m+1}^{J_m})$ and so $f_k(\cup \mathcal I_{m}^{I},\cup \mathcal I_{m}^{J})=\cup_{(I_m,J_m)\in \mathcal I_{m}^{I}\times \mathcal I_{m}^{J}}f_k(\cup \mathcal I_{m+1}^{I_m},\cup \mathcal I_{m+1}^{J_m})=f_k(\cup \mathcal I_{m+1}^{I},\cup \mathcal I_{m+1}^{J})$. Then by Lemma \ref{lemma1} we have
\begin{align*}
  f_k(I,J)&=\cap_{i=n}^\infty f_k(\cup \mathcal I_{i+1}^{I},\cup \mathcal I_{i+1}^{J})\\
  &=f_k(\cap_{i=n}^\infty\cup \mathcal I_{i+1}^{I},\cap_{i=n}^\infty\cup \mathcal I_{i+1}^{J})=f_k(C_\lambda \cap I, C_\lambda \cap J).
\end{align*}

\end{proof}

\begin{lemma}\label{lemma33} Let $ I_1 = J_1=[1-\lambda, 1]$. If $\lambda_k\le \lambda<\frac{1}{2}$, then $f_k(C_\lambda\cap I_1 ,C_\lambda\cap J_1)=[(1-\lambda)^{k+1},1]$.

\end{lemma}
\begin{proof}
Thanks to Lemma \ref{lemma32} it suffices to check that $I_1, J_1$ satisfy \eqref{e33}, that is, we need to verify the following inequalities:
\begin{align*}
    \frac{(1-2\lambda)[1+(k-1)\lambda]}{k\lambda} \le 1-\lambda \qtq{and} 1 \le  \frac{1-\lambda}{k(1-2\lambda)}.
\end{align*}
Note that $\lambda\ge \lambda_k\ge \lambda_2>1/3$, then we have $(k-2)\lambda^2+3\lambda-1\ge0$, this is equivalent to
\begin{equation*}
\frac{(1-2\lambda)[1+(k-1)\lambda]}{k\lambda} \le 1-\lambda.
\end{equation*}
On the other hand, if $\lambda \ge  \frac{k-1}{2k-1}$ we can deduce that $ 1 \le  \frac{1-\lambda}{k(1-2\lambda)}$.

We claim that $\lambda\ge \lambda_k >  \frac{k-1}{2k-1}$. Denote by $x_k=\frac{k-1}{2k-1}$ and $h_k(\lambda)=\frac{(k-1)\lambda^k}{k+1}+2\lambda-1$.  We point that  $h_k(\lambda)$ is increasing, and $h_k(x_k)<0$ can be verified by
the following calculation
\begin{align*}
 h_k(x_k)&=\frac{k-1}{k+1}\left(\frac{k-1}{2k-1}\right)^k+2\left(\frac{k-1}{2k-1}\right)-1\\
         &=\frac{k-1}{k+1}\left(\frac{k-1}{2k-1}\right)^k-\frac{1}{2k-1}\\
         &=\frac{1}{(k+1)(2k-1)}\left[\frac{(k-1)^{k+1}}{(2k-1)^{k-1}}-k-1\right].
\end{align*}
Since $2k-1<2^k<x_k^{-k}$ for all $k\in \nn_{\ge2}$, which implies
\begin{align*}
\frac{(k-1)^{k+1}}{(2k-1)^{k-1}}<k-1<k+1.
\end{align*}
Hence  $ \frac{1}{3}\le x_k<\lambda_k$ and \eqref{e33} holds when $\lambda_k\le \lambda<\frac{1}{2}$.

\end{proof}

We recall our second main result:
\begin{theorem}\label{theorem-22} Given $ k\in \nn_{\ge2}$ and if $\lambda_k\le \lambda<\frac{1}{2}$, then $\set{x^ky:x,y\in C_\lambda}$ is exactly interval $[0,1]$, where $\lambda_k$ is the appropriate root of $(k-1)\lambda^k+(2k+2)\lambda-k-1=0$. Moreover, the sequence $(\lambda_k)$ is increasing to $1/2$ as $k\to \infty$.
\end{theorem}

\begin{proof}

By Lemma \ref{lemma33} $f_k(C_\lambda\cap I ,C_\lambda\cap I)=[(1-\lambda)^{k+1},1]$, where $I=[1-\lambda,1]$. Then $f_k(C_\lambda\cap I ,\lambda(C_\lambda\cap I))=[\lambda(1-\lambda)^{k+1},\lambda]$. Note that $g_k(\lambda)=(1-\lambda)^{k+1}-\lambda$ is decreasing, we denote the appropriate root of $g_k(\lambda)=0$ by $r_k$. This implies that $r_{k+1}<r_k\le r_2(<0.32)<\lambda_2(>0.46)\le \lambda_k$ for any $k\ge2$. We can deduce that $(1-\lambda)^{k+1}<\lambda$ if $\lambda_k\le \lambda<\frac{1}{2}$, and it follows from $C_\lambda=\cup_{n=0}^\infty \lambda^n(C_\lambda\cap I)\cup \set{0}$ that
\begin{equation*}
  \bigcup_{n=0}^\infty f_k(C_\lambda\cap I ,\lambda^n(C_\lambda\cap I))\cup \set{0}=[0,1]\subset f_k(C_\lambda ,C_\lambda)\subset [0,1].
\end{equation*}
Furthermore, we have $h_k(\lambda)>h_{k+1}(\lambda)$, where $h_k(\lambda)$ is defined in the proof of Lemma \ref{lemma33}. We conclude that $(\lambda_k)\nearrow \frac{1}{2} $, and completes the proof.

\end{proof}

\section{Final remarks }\label{S4}
In our work, we first investigate the set $S_t=\set{(x,y)\in C_\lambda\times C_\lambda: xy=t}$, where $t\in(0,1)$. We show that there exists a range of $\lambda$ for which $S_t$ is of cardinality
continuum for any $t\in(0,1)$. The constant in Theorem \ref{theorem-1} admits further improvement. In the first step of the proof of Theorem \ref{theorem-1} we show that $S_t$ is of cardinality
continuum for any $ab<t<(a+\lambda^{n})(b+\lambda^{n})$. In order for $t$ to run over $(0,1)$, we have to choose such basic intervals $ I = [a, a + \lambda^n]  , J = [b, b + \lambda^n]$ satisfy
\begin{equation}\label{p2}
   \frac{(1-\lambda-\lambda^2)(a+\lambda^n)}{\lambda}\le b<b+\lambda^n \le  \frac{a\lambda}{1-2\lambda},
\end{equation}
and $ a + \lambda^n, b + \lambda^n$ are sufficiently close to $1$. So according to \eqref{p2} the constant in Theorem \ref{theorem-1} can be optimized  at most close to $\sqrt{2}-1$ by our method.

This idea can be implemented for other curves. By the same argument as in the proof of Theorem \ref{theorem-1}, we may give similar result for circle:

If $0.459<\lambda<1/2$, then $\set{(x,y):x^2+y^2=r}\cap (C_\lambda\times C_\lambda)$  has cardinality
continuum for any $r\in(0,2)$.
\begin{proof}
Let $D=\set{(x,y):x^2+y^2=r}\cap (C_\lambda\times C_\lambda)$. We recall that \cite[Proposition 2.9]{JKLW}:\\
 $D$  has the cardinality of the
continuum for any $a^2+b^2<r<(a+\lambda^n)^2+(b+\lambda^n)^2$, where $I=[a,a+\lambda^n],J=[b,b+\lambda^n]$ with $a\le b$ satisfy \eqref{p2} (by the same argument as in the Case (B)(ii) of proof of Theorem \ref{theorem-1}, we can deduce that $a$ and $b$ may be equal).

Next we choose $I_{2,1}=J_{2,1}=[1-\lambda^2,1];I_{2,2}=[1-\lambda,1-\lambda+\lambda^2],J_{2,2}=[1-\lambda^2,1];I_{2,3}=J_{2,3}=[1-\lambda,1-\lambda+\lambda^2]; I_{2,4}=[\lambda-\lambda^2,\lambda],J_{2,4}=[1-\lambda,1-\lambda+\lambda^2]$.
$I_{2,1},J_{2,1}$ satisfy \eqref{p2} for $\lambda>0.446$, $I_{2,2},J_{2,2}$ satisfy \eqref{p2} for $\lambda>0.400$, $I_{2,3},J_{2,3}$ satisfy \eqref{p2} for $\lambda>0.459$ and $I_{2,4},J_{2,4}$ satisfy \eqref{p2} for $\lambda>0.431$.
On the other hand, $(1-\lambda+\lambda^2)^2+1>2(1-\lambda^2)^2$ for $\lambda>0.320$, $2(1-\lambda+\lambda^2)^2>(1-\lambda)^2+(1-\lambda^2)^2$ for $\lambda>0.350$, $\lambda^2+(1-\lambda+\lambda^2)^2>2(1-\lambda)^2$ for $\lambda>0.394$, by \cite[Proposition 2.9]{JKLW} $D$ has the cardinality of the
continuum for any $(\lambda-\lambda^2)^2+(1-\lambda)^2<r<2$. Furthermore,  $(\lambda-\lambda^2)^2+(1-\lambda)^2<2\lambda^2$ for $\lambda>0.436$, so if $0.459<\lambda<1/2$ it follows that $D$ is of cardinality
continuum for any $r\in(0,2)$.

\end{proof}

In the second part of paper, we consider the image of $C_\lambda\times C_\lambda$ under the continuous function $f_k(x,y)=x^ky$, $k\in\nn_{\ge2}$. Note that $f_1(C_{1/3}\times C_{1/3})\subsetneq[0,1]$, $f_2(C_{1/3}\times C_{1/3})=[0,1]$(\cite{ART,JK}). What about the case for $k\ge3$, do we have $f_k(C_{1/3}\times C_{1/3})=[0,1]$ for all $k\ge2$?

For $C_\lambda$ the basic intervals of rank $2$ are
\begin{equation*}
  [0,\lambda^2], [\lambda-\lambda^2,\lambda], [1-\lambda,1-\lambda+\lambda^2], [1-\lambda^2,1].
\end{equation*}
Denote by $A$ the union of these intervals. It follows that
\begin{align*}
  f_k(A,A)&=
 [0,\lambda^{2}]\cup  [(\lambda-\lambda^2)(1-\lambda)^k, \lambda (1-\lambda+\lambda^2)^k]\\
&  \cup [(1-\lambda)^{k+1}, (1-\lambda+\lambda^2)^{k+1})] \cup  [(1-\lambda^2)(1-\lambda)^k, (1-\lambda+\lambda^2)^k]\\
&\cup [(\lambda-\lambda^2)(1-\lambda^2)^k, \lambda]\cup [(1-\lambda)(1-\lambda^2)^k, 1-\lambda+\lambda^2]\\
 & \cup [(1-\lambda^2)^{k+1}, 1].
 \end{align*}
If $\lambda<(1-\lambda+\lambda^2)^k<(1-\lambda)(1-\lambda^2)^k$, we see $((1-\lambda+\lambda^2)^k,(1-\lambda)(1-\lambda^2)^k)$ is missing from $[0,1]$.  For example,
when $\lambda=1/3$: $f_3(A,A)=[0,1]$; and $f_4(A,A)=[0,\frac{2401}{19683}]\cup  [ \frac{32}{243},\frac{2401}{6561}] \cup [\frac{8192}{19683},1]$.
The intervals $(\frac{2401}{19683},\frac{32}{243})$ and $(\frac{2401}{6561},\frac{8192}{19683})$ are missing from $[0,1]$. So $f_4(C_{1/3}\times C_{1/3})\subsetneq[0,1]$. Our observations reveal that as
$k(\ge2)$ increases, the set of $\lambda$ for $f_k(C_{\lambda}\times C_{\lambda})=[0,1]$ becomes smaller.

At the end, we are interested in the question: the constants in Theorems \ref{theorem-1} and \ref{theorem-2} are not optimal. So can we determine the critical values in Theorems \ref{theorem-1} and \ref{theorem-2}?


\vspace{1cm}




{\bf Acknowledgments}

The authors are grateful to Zhiqiang Wang for valuable discussions on our paper. Y. Cai was supported by National Natural Science Foundation of China No. 12301108. L. Wang was supported by Haizhou Bay Talents Innovation Program No. KQ25011.



 \end{document}